\documentclass[12pt]{article}

\usepackage{amsmath, epsfig, cite, lineno}
\usepackage{amssymb,amsthm}
\usepackage{amsfonts, color}
\usepackage{latexsym}

\newtheorem{thm}{Theorem}

\newtheorem{conj}[thm]{Conjecture}
\newtheorem{lem}[thm]{Lemma}

\newcommand{\pf}{\noindent{\it Proof.} }

\numberwithin{equation}{section}

\begin{document}


\begin{center}
{\Large\bf Proof of a conjecture of Mircea Merca}
\end{center}

\vskip 2mm \centerline{Victor J. W. Guo}
\begin{center}
{\footnotesize Department of Mathematics, Shanghai Key Laboratory of
PMMP, East China Normal University,\\ 500 Dongchuan Rd., Shanghai 200241,
 People's Republic of China\\
{\tt jwguo@math.ecnu.edu.cn,\quad http://math.ecnu.edu.cn/\textasciitilde{jwguo}}}
\end{center}


\vskip 0.7cm \noindent{\bf Abstract.} We prove that, for any prime $p$ and positive integer $r$ with $p^r>2$,
the number of multinomial coefficients such that
$$
{k\choose k_1,k_2,\ldots,k_n}=p^r,\quad \text{and}\quad k_1+2k_2+\cdots+nk_n=n,
$$
is given by
$$
\delta_{p^r,\,k}\left(\left\lfloor\frac{n-1}{p^r-1}\right\rfloor -\delta_{0,\,n\bmod p^r} \right),
$$
where $\delta_{i,j}$ is the Kronecker delta and $\lfloor x\rfloor$ stands for the largest integer not exceeding $x$.
This confirms a recent conjecture of Mircea Merca.



\section{Introduction}
The multinomial coefficients are
defined by
$$
{k\choose k_1,k_2,\ldots,k_n}=\frac{k!}{k_1!k_2!\cdots k_n!},
$$
where $k=k_1+k_2+\cdots+k_n$. Fine \cite[p.~87]{Fine} gave a connection between binomial coefficients and binomial coefficients:
\begin{align}
\sum_{\substack{k_1+k_2+\cdots+k_n=k\\ k_1+2k_2+\cdots+nk_n=n}}{k\choose k_1,k_2,\ldots,k_n}={n-1\choose k-1}. \label{eq:fine}
\end{align}
Let $M_m(n,k)$ be the number of multinomial coefficients such that
$$
{k\choose k_1,k_2,\ldots,k_n}=m,\quad \text{and}\quad k_1+2k_2+\cdots+nk_n=n.
$$
For example, we have $M_{6}(10,3)=4$, since
$$
10=1+2+7=1+3+6=1+4+5=2+3+5.
$$
It is easy to see that $M_{1}(n,k)=\delta_{0,\,n\bmod k}$.
Recently, applying Fine's formula \eqref{eq:fine}, Merca \cite{Merca} obtained new upper bounds involving $M_m(n,k)$ for the number of partitions of $n$ into $k$ parts.
He also proved that
\begin{align*}
M_{2}(n,k)&=\delta_{2,\,k}\left\lfloor\frac{n-1}{2}\right\rfloor, \\
M_{p}(n,k)&=\delta_{p,\,k}\left(\left\lfloor\frac{n-1}{p-1}\right\rfloor -\delta_{0,\,n\bmod p}\right),
\end{align*}
where $p$ is an odd prime.

In this paper, we shall prove the following result, which was conjectured by Merca \cite[Conjecture 1]{Merca}.
\begin{thm}\label{thm:1} Let $p$ be a prime and let $n,k,r$ be positive integers with $p^r>2$. Then
$$
M_{p^r}(n,k)=\delta_{p^r,\,k}\left(\left\lfloor\frac{n-1}{p^r-1}\right\rfloor -\delta_{0,\,n\bmod p^r} \right).
$$
\end{thm}
Merca \cite{Merca} pointed out that, when $m$ is not a prime power, the formula for $M_{m}(n,k)$
is more involved. For example, we have
$$
M_{10}(n,k)=\delta_{10,\,k}\left(\left\lfloor\frac{n-1}{9}\right\rfloor-\delta_{0,\,n\bmod 10}\right)
+\delta_{5,\,k}\left(\left\lfloor\frac{n+1}{6}\right\rfloor-\delta_{0,\,n\bmod 5}-\delta_{0,\,n\bmod 6}\right).
$$


\section{Proof of Theorem \ref{thm:1} }
We need the following result.
\begin{lem}\label{lem:1}Let $n$ and $k$ be two positive integers with $2\leqslant k\leqslant \frac{n}{2}$. Then
the binomial coefficient ${n\choose k}$ is not a prime power.
\end{lem}
\pf For any prime $p$, the $p$-adic order of $n!$ can be given by
\begin{align*}
{\rm ord}_p n!=\sum_{i=1}^\infty\left\lfloor\frac{n}{p^i}\right\rfloor.
\end{align*}
If ${n\choose k}=\frac{n!}{k!(n-k)!}$ were a prime power, say $p^r$, then
\begin{align}
r=\sum_{i=1}^\infty\left(\left\lfloor\frac{n}{p^i}\right\rfloor-
\left\lfloor\frac{n-k}{p^i}\right\rfloor-\left\lfloor\frac{n}{p^i}\right\rfloor\right).  \label{eq:r}
\end{align}
Note that $\lfloor x+y\rfloor -\lfloor x\rfloor -\lfloor y\rfloor\leqslant 1$. From \eqref{eq:r} we deduce that $r$ is less than or equal to the largest integer $i$ such that
$p^{i}\leqslant n$. Namely, $p^r\leqslant n$. On the other hand, for $2\leqslant k\leqslant n-2$, we have
${n\choose k}>n$, a contradiction.  Therefore, the initial assumption must be false.
\qed
\medskip

\noindent{\it Proof of Theorem {\rm\ref{thm:1}.}}
Let
\begin{align}
{k_1+k_2+\cdots+k_n\choose k_1,k_2,\ldots,k_n}=p^r.  \label{eq:multi-p}
\end{align}
We assert that there are exactly two $i$'s such that $k_i\geqslant 1$. In fact, if $k_1,k_2,k_3\geqslant 1$,
then either ${k_1+k_2+k_3\choose k_1,k_2,k_3}={3\choose 1,1,1}=6$, or by Lemma \ref{lem:1}, ${k_1+k_2+k_3\choose k_a}$ ($k_a=\max\{k_1,k_2,k_3\}$) is not a prime power.
But this is impossible, since both ${k_1+k_2+k_3\choose k_1,k_2,k_3}$ and ${k_1+k_2+k_3\choose k_a}$
divide ${k_1+k_2+\cdots+k_n\choose k_1,k_2,\ldots,k_n}$. This proves the assertion. Furthermore, by Lemma \ref{lem:1} again,
one of the two non-zero $k_i$'s must be $1$, and by \eqref{eq:multi-p}, the other non-zero term is equal to $p^r-1$.
In other words, the identity \eqref{eq:multi-p} holds if and only if $(k_1,k_2,\ldots,k_n)$
is a rearrangement of $(p^r-1,1,0,\ldots,0)$. Consider the equation
\begin{align}
(p^r-1)x+y=n.  \label{eq:xy}
\end{align}
If $k=p^r$, then we conclude that $M_{p^r}(n,k)$ is equal to the number of solutions $(x,y)$ to \eqref{eq:xy} with $x\neq y$,
i.e.,
$$
M_{p^r}(n,k)=\left\lfloor\frac{n-1}{p^r-1}\right\rfloor
-\begin{cases}
1, &\text{if $n\equiv 0\pmod{p^r}$,}\\
0, &\text{otherwise.}
\end{cases}
$$
If $k\neq p^r$, then it is obvious that $M_{p^r}(n,k)=0$. This completes the proof. \qed

\section{Concluding remarks}
Note that, Lemma \ref{lem:1} plays an important part in our proof of Theorem \ref{thm:1}.
It seems that we may say something more about the factors of ${n\choose k}$ for $2\leqslant k\leqslant \frac{n}{2}$.
Since ${n\choose k}=\frac{n}{k}{n-1\choose k-1}$, we have $\gcd\left({n\choose k},n\right)>1$, where $\gcd(a,b)$
denotes the greatest common divisor of two integers $a$ and $b$. If
\begin{align}
\gcd\left({n\choose k},n-1\right)>1, \label{eq:nk-1}
\end{align}
then, noticing that $\gcd(n,n-1)=1$, we immediately deduce that ${n\choose k}$ has at least two different
prime factors (namely, Lemma \ref{lem:1} holds). But, in general, the inequality \eqref{eq:nk-1} does
not hold. For example, $\gcd\left({7\choose 3},6\right)=1$. Similarly, the identity $\gcd\left({14\choose 4},12\right)=1$ means that
we cannot expect
\begin{align}
\gcd\left({n\choose k},n-2\right)>1. \label{eq:nk-2}
\end{align}
We close our paper with the following conjecture, which asserts that at least one of \eqref{eq:nk-1} and \eqref{eq:nk-2}
is true.
\begin{conj}
Let $n$ and $k$ be two positive integers with $2\leqslant k\leqslant \frac{n}{2}$. Then
$$
\gcd\left({n\choose k},{n-1\choose 2}\right)>1.
$$
\end{conj}
We have verified the above conjecture for $n$ up to $5000$ via Maple.


\vskip 5mm
\noindent{\bf Acknowledgments.} This work was partially
supported by the Fundamental Research Funds for the Central Universities and
the National Natural Science Foundation of China (grant 11371144).

\end{document}